\documentclass[letterpaper, 10 pt, conference]{ieeeconf}  

\IEEEoverridecommandlockouts                              
\title{\LARGE \bf	
Robust synchronization of heterogeneous robot swarms on the sphere
}

\author{Johan Markdahl, Daniele Proverbio, and Jorge Goncalves
\thanks{This work is part of the Luxembourg National Research Fund (FNR) project SYBION. The authors gratefully acknowledge having received financial support from the FNR through their CORE OPEN instrument.}
\thanks{The authors are with the Luxembourg Centre for Systems Biomedicine at the University of Luxembourg. Email: {\tt\small markdahl@kth.se}}%
}

\usepackage{graphicx}

\usepackage{subdepth} 

\usepackage{graphics} 
\usepackage{epsfig} 

\usepackage{amssymb} 
\usepackage{amsmath,psfrag,color,graphicx} 
\usepackage{mathabx}

\usepackage{thmtools}
\renewcommand\thmcontinues[1]{Continued}

\usepackage{graphicx}
\usepackage{quoting} 
\usepackage{lettrine} 
\usepackage{tocloft} 
\usepackage{titletoc}
\usepackage{fmtcount}
\usepackage{calc} 
\usepackage{cite} 
\usepackage{multicol} 
\usepackage{centernot}
\usepackage{nicefrac}
\usepackage{dsfont}
\usepackage{booktabs}
\usepackage{mathtools, cuted}

\usepackage{lipsum}
\usepackage{psfrag}
\usepackage[export]{adjustbox} 

\usepackage{mathtools}

\PassOptionsToPackage{cmyk}{xcolor}
\usepackage[cmyk]{xcolor}
\selectcolormodel{cmyk}

\usepackage{textcomp}

\usepackage{txfonts}
\usepackage{lmodern}
\usepackage{times}

\DeclareFontEncoding{LS1}{}{}
\DeclareFontSubstitution{LS1}{stix}{m}{n}

\makeatletter
\newcommand*\textmathversion{\csname textmv@\math@version\endcsname}
\newcommand*\textmv@normal{m}
\newcommand*\textmv@bold{b}
\makeatother

\newcommand{\ve}[2][]{\ensuremath{\boldsymbol{\mathrm{#2}}}_{#1}}
\newcommand{\vet}[2][]{\ensuremath{\smash{\boldsymbol{\mathrm{#2}}^{\!\top}_{#1}}}}
\newcommand{\vd}[2][]{\ensuremath{\dot{\boldsymbol{\mathrm{#2}}}_{#1}}}

\newcommand{\ma}[2][]{\ensuremath{\boldsymbol{\mathrm{#2}}}_{#1}}
\newcommand{\mat}[2][]{\ensuremath{\boldsymbol{\mathrm{#2}}^{\!\top}_{#1}}}

\DeclareMathOperator{\grad}{\ensuremath{\mathrm{grad}}}

\newcommand{\R}{\ensuremath{\mathds{R}}}

\newcommand{\C}{\ensuremath{\mathds{C}}}

\newcommand{\N}{\ensuremath{\mathds{N}}}

\newcommand{\so}{\ensuremath{\mathsf{so}(n)}}

\newcommand{\ie}{\textit{i.e.}, }

\newcommand{\eg}{\textit{e.g.}, }

\newcommand{\mtr}{\hspace{-0.3mm}\ensuremath{^\top}}

\newcommand{\qedsymboltext}{\hspace{\fill}$\blacksquare$} 

\newcommand\raiseT[2]{\setbox0\hbox{$#1{#2}$}\raise\dp0\box0}

\newcommand{\V}{\ensuremath{\mathcal{V}}}

\newcommand{\E}{\ensuremath{\mathcal{E}}}

\newcommand{\etc}{\emph{etc.\ }}

\DeclareMathOperator{\diag}{\ensuremath{\mathrm{diag}}}

\DeclareMathOperator{\trace}{\ensuremath{\mathrm{tr}}}

\DeclareMathOperator{\im}{\ensuremath{\mathrm{Im}}}

\DeclareMathOperator{\re}{\ensuremath{\mathrm{Re}}}

\newcommand{\Ni}{\ensuremath{\mathcal{N}_i}}

\makeatletter
\newcommand{\raisemath}[1]{\mathpalette{\raisem@th{#1}}}
\newcommand{\raisem@th}[3]{\raisebox{#1}{$#2#3$}}
\makeatother

\newcommand{\ts}[2][]{\ensuremath{\mathsf{T}_{#2}#1}}

\definecolor{kthbluergb}{RGB}{25,84,166}
\definecolor{kthbluecmyk}{cmyk}{1,0.55,0,0}

\definecolor{kthblueA}{RGB}{25,84,166}
\definecolor{kthblueB}{RGB}{46,124,192}
\definecolor{kthblueC}{RGB}{112,153,209}
\definecolor{kthblueD}{RGB}{164,186,225}
\definecolor{kthblueE}{RGB}{211,220,241}

\newcounter{counter} 

\newtheorem{theorem}[counter]{Theorem}

\newtheorem{definition}[counter]{Definition}

\newtheorem{remark}[counter]{Remark}


\newcounter{parentnumber}

\begin{document}

\maketitle
\thispagestyle{empty}
\pagestyle{empty}

\begin{abstract}
Synchronization on the sphere is important to certain control applications in swarm robotics. Of recent interest is the Lohe model, which  generalizes the Kuramoto model from the circle to the sphere. The Lohe model is mainly studied in mathematical physics as a toy model of quantum synchronization. The model makes few assumptions, wherefore it is well-suited to represent a swarm. Previous work on this model has focused on the cases of complete and acyclic networks or the homogeneous case where all oscillator frequencies are equal. This paper concerns the case of heterogeneous oscillators connected by a non-trivial network. We show that any undesired equilibrium is exponentially unstable if the frequencies satisfy a given bound. This property can also be interpreted as a robustness result for  small model perturbations of the homogeneous case with zero frequencies. As such, the Lohe model is a good choice for control applications in swarm robotics.
	



\end{abstract}

\section{INTRODUCTION}

\noindent Synchronization on nonlinear spaces is important to a number of robotics applications. Many such problem pose special challanges due to the non-Euclidean topology of the state space, which also make them particularly interesting from a control theory perspective \cite{bhat2000topological,sepulchre2011consensus,markdahl2021tac}. Consensus protocols have been developed to achieve synchronization on general manifolds \cite{sarlette2009consensus,tron2013riemannian}. However, for particular manifolds, other algorithms can outperform those protocols \cite{markdahl2018tac,markdahl2020high}. This includes the so-called Lohe model which is studied in  physics \cite{lohe2009non,lohe2010quantum,zhu2013synchronization,chi2014emergent,crnkic2018swarms,crnkic2019data,chandra2019continuous,chandra2019complexity}. The Lohe model is a generalization of the Kuramoto model of a multi-agent system of coupled oscillators from the circle to the sphere. Compared to other methods \cite{sarlette2009consensus}, the Lohe model is less demanding in terms of computation and communication. As such, it is suitable for robot swarms where each robot has limited resources available for control. We have previously proved that the Lohe model displays almost global synchronization \cite{markdahl2018tac,markdahl2020high}. In this paper we show that this property is also robust under small drift-like  perturbations. Our main result amounts to a condition on the model parameters under which undesired equilibrium points remain exponentially unstable.


Roughly speaking, our result can be interpreted as robustness of almost global synchronization. The homogeneous Lohe model is a driftless system. It becomes a system with drift (the heterogenous model) after a small perturbation term is added to the dynamics. Almost global synchronization becomes almost global practical synchronization \cite{montenbruck2015synchronization} (up to a technicality). The property of almost global sync is important for two main reasons. First, in reduced rigid-body attitude synchronization, the guarantee of almost global sync is clearly preferable to local sync. In particular, conditions for local convergence to synchronization on $\mathcal{S}^n$ requires all initial states to belong to an open hemisphere \cite{zhu2013synchronization,thunberg2018lifting}. The probability of drawing such a configuration from a uniform distribution of all configurations on $(\mathcal{S}^n)^N$ goes to zero exponentially as the number of agents $N$ goes to infinity. By contrast, almost global  sync does not depend on $N$. For almost global practical sync, we again find a dependence on $N$ in our result. Secondly, almost global practical synchronization for all connected networks distinguishes the heterogenous Lohe model from the heterogenous Kuramoto, model which is multistable. The fact that the two models behave qualitatively different is intruiging and motivates further study of the heterogenous Lohe model.

\subsection{Literature review}

\noindent The Kuramoto model is a popular model of collective behavior in multi-agent systems of weakly coupled oscillators \cite{hoppensteadt2012weakly,dorfler2014synchronization,rodrigues2016kuramoto}. Recent years have seen a growing interest in high-dimensional generalizations of the Kuramoto model, the most popular of which is the Lohe model of synchronization on $\mathcal{S}^n$ \cite{olfati2006swarms,lohe2009non,lohe2010quantum,zhu2013synchronization,chi2014emergent,markdahl2018tac,crnkic2018swarms,crnkic2019data,chandra2019continuous,chandra2019complexity,markdahl2020high}. The study of this model in physics is  motivated by its relation to synchronization of quantum bits \cite{lohe2009non,lohe2010quantum}, but it also appears in control systems designed to achieve reduced rigid-body attitude synchronization, \ie to coordinate the pointing orientations of robots \cite{song2017intrinsic,markdahl2018tac,casau2018}, in bio-inspired models of source-seeking and learning \cite{al2018gradient,crnkic2018swarms}, and in machine learning applications \cite{crnkic2019data}. There are several variations of the model, including second-order dynamics \cite{li2014unified,kim2020asymptotic} and discrete-time maps \cite{li2015collective}. A limitation in our current understanding of the Lohe model is the restriction to combinations of complete or acyclic networks \cite{chi2014emergent,casau2018,chandra2019continuous,chandra2019complexity,kim2020asymptotic}, homogeneous frequencies \cite{olfati2006swarms,zhu2013synchronization,li2014unified,markdahl2018tac}, and local behavior \cite{zhu2013synchronization,li2014unified,chi2014emergent,lageman2016consensus,kim2020asymptotic}. By contrast, this paper concerns the  global behavior of the Lohe model with heterogeneous frequencies over non-trivial networks. Previous work established almost global sync over all connected networks for the homogeneous Lohe model on the $n$-sphere for $n\geq2$ \cite{markdahl2018tac}. This paper generalizes the results of  \cite{markdahl2018tac} to the case of heterogeneous frequencies under a condition on the model parameters.



The literature contains a result for the homogeneous Lohe model under small homogenous and heterogeneous perturbations \cite{lageman2016consensus}. However, it is limited to local existence of asymptotically stable sets without further characterization. This paper studies the effect of perturbations on a global level and gives an explicit bound on a range of parameter values that excludes the possibility of critical transitions from instability to stability. There is also a robust hybrid feedback algorithm for global synchronization on the sphere \cite{casau2018}, but convergence results  are limited to tree graphs. Furthermore, there are many general results about small perturbations in the control theory literature \cite{khalil2002nonlinear}, including applications to synchronization of heterogeneous agents \cite{montenbruck2015synchronization}. Moreover, exponentially stable systems are known to be robust in general  \cite{khalil2002nonlinear}. In terms of perturbation theory, the contribution of this paper is to provide an explicit relation between the parameters of the system that, if satisfied, guarantees exponential instability of undesired equilibria.





 
 

\section{PRELIMINARIES}

\noindent The following notation is used in this paper. The inner product of $\ve{x}, \ve{y}\in\R^{n}$ is denoted by $\langle\ve{x},\ve{y}\rangle$. Let $\|\cdot\|_2$ denote the Euclidean norm of a vector or induced $2$-norm of a matrix. The $n$-sphere is $\mathcal{S}^n=\{\ve{x}\in\R^{n+1}\,|\,\|\ve{x}\|_2=1\}$. The tangent space of $\mathcal{S}^n$ is  $\ts[\mathcal{S}^n]{\ve{x}}=\{\ve{y}\in\R^{n+1}\,|\,\langle\ve{y},\ve{x}\rangle=0\}\simeq\R^{n}$. The special orthogonal lie algebra is $\so=\{\ma{S}\in\R^{n+t\times n+1}\,|\,\mat{S}=\ma{S}\}$. The gradient in Euclidean space is denoted $\nabla:f(\ve{x})\mapsto\nabla f(\ve{x})\in\R^{n+1}$, the intrinsic gradient map at a point $\ve{x}\in\mathcal{S}^n$ is denoted $\grad:f(\ve{x})\mapsto\grad f(\ve{x})\in\ts[\mathcal{S}^n]{\ve{x}}$.  

An undirected, simple \emph{graph} is a pair $\mathcal{G}=(\mathcal{V},\mathcal{E})$ where $\mathcal{V}\subset\N$ is the node set and $\mathcal{E}\subset \{e\subset\V\,|\,|e|=2\}$ is the edge set. The graph $\mathcal{G}$ encodes the communication topology of a multi-agent system where each $i\in\V$ corresponds to an agent and each $\{i,j\}\in\E$ indicates a bidirectional information flow between agent $i$ and agent $j$. Throughout the paper it is assumed that $\mathcal{G}$ is connected. The set of neighbors of agent $i$ is $\mathcal{N}_i=\{j\in\V\,|\,\{i,j\}\in\E\}$. Other objects associated with agent $i$ also carry $i$ as a subindex, \eg the state $\ve[i]{x}\in\mathcal{S}^n$, the tangent space of $\mathcal{S}^n$ at $\ve[i]{x}$ is $\ts[\mathcal{S}^n]{i}:=\ts[\mathcal{S}^n]{\ve[i]{x}}$ \etc

The exponential instability property of an equilibrium is important for stability analysis. It appears in the the indirect method of Lyapunov, the Hartman-Grobman theorem, and the center manifold theorem \cite{sastry2013nonlinear}. In the proof of our main result we use the center manifold theorem to show that exponential instability, defined in terms of the spectral abscissa of the linearization matrix, relates to a lack of Lyapunov stability. For now, we just give a formal definition:

\begin{definition}
Let $\ma{A}(\ve{x})$ denote the linearization matrix of a continuous-time, autonomous dynamical system around an equilibrium $\ve{x}$. The \emph{spectral abscissa} of $\ma{A}(\ve{x})$ is
\begin{align*}
\eta(\ve{x})=\!\!\!\!\!\!\max_{\alpha(\ve{x})\in\sigma(\ma{A}(\ve{x}))}\!\!\!\!\!\!\re\alpha(\ve{x}),
\end{align*}
where $\sigma(\cdot)$ denotes the spectrum of a matrix. The equilibrium $\ve{x}$ is \emph{exponentially unstable} if $\eta(\ve{x})>0$.
\end{definition}

\subsection{The Lohe model}

\noindent The \emph{heterogeneous Lohe model on networks} is given by
\begin{align}\label{eq:dxidt}
\vd[i]{x}&=\ma[i]{\Omega}\ve[i]{x}+(\ma[n]{I}-\ve[i]{x}\vet[i]{x})\sum_{j\in\Ni}k_{ij}\ve[j]{x},
\end{align}
where $\ma[i]{\Omega}\in\mathsf{so}(n)$, $\ve[i]{x}\in\mathcal{S}^n$, and $k_{ij}\in(0,\infty)$, $k_{ij}=a_{ji}$ are the coupling gains. The $n$-sphere is invariant under these dynamics by Nagumo's invariance theorem since the right-hand side of \eqref{eq:dxidt} belongs to the tangent space $\ts[\mathcal{S}^n]{i}$ \cite{blanchini2008set}. For $n=1$ we obtain the \emph{Kuramoto model on networks} \cite{rodrigues2016kuramoto} after a change from Cartesian to polar coordinates,
\begin{align}\label{eq:kuramoto}
\dot{\theta}_i=\omega_i+\sum_{j\in\Ni}k_{ij}\sin(\theta_j-\theta_i),
\end{align}
where $\cos\theta_i=\langle\ve[i]{x},\ve[1]{e}\rangle$, $\omega_i=\langle\ve[1]{e},\ma[i]{\Omega}\ve[2]{e}\rangle$. 

The \emph{homogeneous Lohe model on networks} is obtained by setting $\ma[i]{\Omega}=\ma{\Omega}$ and is equivalent to the case of $\ma[i]{\Omega}=\ma{0}$ \cite{markdahl2020high},
\begin{align}\label{eq:homo} 
\vd[i]{z}=(\ma{I}-\ve[i]{z}\vet[i]{z})\sum_{j\in\Ni}k_{ij}\ve[j]{z}.
\end{align}
The model \eqref{eq:homo} is the gradient flow of the disagreement  function $V:(\mathcal{S}^n)^N\rightarrow[0,\infty)$ which is given by
\begin{align*}
V:=\tfrac12\sum_{i\in\V}\sum_{j\in\Ni}k_{ij}\|\ve[i]{z}-\ve[j]{z}\|^2_2.
\end{align*}
Denote $\ve{z}:=(\ve[i]{z})_{i=1}^N$. Then 
\begin{align}\label{eq:homo2}
\vd{z}=-\grad V(\ve{z})=-(\grad_iV(\ve{z}))_{i=1}^N,
\end{align}
where $\grad$ and $\grad_i$ denotes the gradient on the sphere with respect to $\ve{z}$ and $\ve[i]{z}$ respectively. The gradient $\grad_i V(\ve{z})$ is calculated in \cite{markdahl2018tac} by taking the gradient with respect to $\ve[i]{z}\in\R^{n+1}$, $\nabla_i V=\sum_{j\in\Ni}k_{ij}\ve[j]{z}$, and applying an orthogonal projection operator $\ma[i]{P}:\R^{n+1}\rightarrow\ts[\mathcal{S}^n]{i}$ to it,
\begin{align*}
\grad_iV(\ve{z})=\ma[i]{P}\nabla_i V(\ve{z})=(\ma{I}-\ve[i]{z}\vet[i]{z})\sum_{j\in\Ni}k_{ij}\ve[j]{z}.
\end{align*}
%


The heterogeneous Lohe model \eqref{eq:dxidt} can be written as the sum of a \emph{drift} term and a \emph{gradient descent flow} term,
\begin{align}\label{eq:flow}
\vd[i]{x}&=\ma[i]{\Omega}\ve[i]{x}-\grad_i V(\ve{x}),
\end{align}
The drift term $\ma[i]{\Omega}\ve[i]{x}$ cannot be written as the gradient of a scalar field, since it would result in the Hessian $\ma[i]{\Omega}$ being skew-symmetric. This, in turn, would contradict Schwarz's theorem on the equality of mixed partial derivatives. For the Kuramoto model \eqref{eq:kuramoto} it is however possible to write this term as the gradient of a potential function, see \cite{mirollo2005spectrum}.

From a control theory perspective, the model \eqref{eq:dxidt} may be taken to represent a system with input $\ve[i]{u}:\R^{n+1}\rightarrow\ts[\mathcal{S}^n]{i}$, 
\begin{align*}
\vd[i]{x}=\ma[i]{\Omega}\ve[i]{x}+\ve[i]{u},
\end{align*}
where $\ve[i]{x}\in\mathcal{S}^n$, $\ma[i]{\Omega}\in\so$, and the feedback law $\ve[i]{u}=-\grad_iV(\ve{x})$ results in \eqref{eq:flow}. The drift term $\ma[i]{\Omega}\ve[i]{x}$ is assumed to be unknown, but small. As such, it can be interpreted as a form of bias or model error. An observer-based feedback may be able to estimate and cancel the drift term. Here we assume that such a feedback, or some other advanced algorithm (\emph{cf.} \cite{sarlette2009consensus,thunberg2017dynamic}), is not used due to a requirement that robots in the swarm are relatively cheap in terms of sensors, communication ability, and computational power.

\subsection{Notions of synchronization}

\noindent To explain the goal of this paper we need some formal notions of synchronization for the Lohe model.


\begin{definition}[Phase synchronization]
The system \eqref{eq:dxidt} is \emph{phase synchronized} if  $\|\ve[i]{x}-\ve[j]{x}\|_2=0$ for all $\{i,j\}\in\E$.
\end{definition}

Note that if the system is phase synchronized then $\vd[i]{x}=\ma[i]{\Omega}\ve[i]{x}$ whereby the system will break away from phase synchronization unless $\ma[i]{\Omega}=\ma{\Omega}$, \ie unless the system is homogeneous. Since this paper studies the heterogeneous Lohe model we need a different notion of synchronization:


\begin{definition}[Practical synchronization]\label{def:practical}
The system \eqref{eq:dxidt} is \emph{practically synchronized} if all agents are contained in an open convex cone with angle less than $\pi/2$.
\end{definition}


For the homogeneous Kuramoto model there exists asymptotically stable equilibria that are not practically synchronized, \eg the $q$-twisted states
\begin{align*}
\mathcal{S}_q&=\{(\theta_i)_{i=1}^N\in[0,2\pi)^N\,|\,\theta_i=\varphi+\tfrac{2\pi q}{N},\,\varphi\in[0,2\pi)\}.
\end{align*}
By contrast, for the homogeneous Lohe model, the only asymptotically stable equilibrium set is phase synchronization \cite{markdahl2018tac}. In this paper we show that this results extends to the heterogeneous model: for sufficiently small frequencies $\ma[i]{\Omega}$, the only stable configurations is practical synchronization.

Finally, we give a definition of dispersed equilibria. For this we also need the concept of an open hemisphere. The intersection of $\mathcal{S}^n$ and an open halfspace,
\begin{align*}
\mathcal{S}^n\cap\{\ve{x}\in\R^{n+1}\,|\,\langle\ve{x},\ve{y}\rangle>0\}
&=\{\ve{x}\in\mathcal{S}^n\,|\,\langle\ve{x},\ve{y}\rangle>0\},
\end{align*}
is an \emph{open hemisphere} for any $\ve{y}\neq\ve{0}$.

\begin{definition}\label{def:incohesive}
An equilibrium $(\ve[i]{x})_{i=1}^N$ is referred to as \emph{dispersed} if it there is no open hemisphere containing the set $\{\ve[1]{x},\ldots,\ve[N]{x}\}$ of all agents.
\end{definition}

Note that a open hemispheres are the largest \emph{geodesically convex} subset of a sphere, \ie there is a unique geodesic between any pair of points in the set. As such, open hemispheres are \emph{homotopic} (topologically equivalent) to $\R^n$. There are a number of result about synchronization of various Lohe models on hemispheres \cite{zhu2013synchronization,thunberg2018lifting}. However, due to this topological difference, the global analysis of dispersed equilibria presents additional challanges \cite{bhat2000topological}. 

\section{MAIN RESULTS}

\noindent Because the full proof of our main result is long, it is easy to get lost in the details. We provide a brief proof sketch to explain the main ideas. All the details are in the Appendix.
%
%


%
%

\begin{theorem}\label{th:main}
Any dispersed equilibrium of the heterogeneous Lohe model \eqref{eq:dxidt} on $\mathcal{S}^n$ is exponentially unstable if the frequencies $\ma[i]{\Omega}$ are small in the following sense
\begin{align}\label{eq:th}
\bigl(\sum_{i\in\V}\!\|\ma[i]{\Omega}\|^2_2\bigr)^{\frac12}<\tfrac{K}{n+1}(n-1-\cos\tfrac{\pi}{N})(1-\cos\tfrac{\pi}N),
\end{align}
where $K=\min_{\{i,j\}\in\E} k_{ij}$ and $n\geq2$. The bound \eqref{eq:th} is in $\mathcal{O}(K/N^4)$ for $n=2$ and $\mathcal{O}(K/N^2)$ for $n\geq3$. 
\end{theorem}

\emph{Proof sketch:} \noindent Let $\ma{A}(\ve{x})$ denote the linearization of \eqref{eq:dxidt} at a dispersed equilibrium $\ve{x}\in\mathcal{S}^n$. The center manifold theorem can be used to establish the instability of $\ve{x}$ when the \emph{spectral abscissa}  of $\ma{A}(\ve{x})$, \ie the largest real part $\re\alpha(\ve{x})$ of all eigenvalues, is positive \cite{sastry2013nonlinear}. The linearization matrix $\ma{A}(\ve{x})$ is the sum of a skew-symmetric and a symmetric  matrix,
\begin{align*}
\ma{A}(\ve{x})=\diag(\ma[1]{\Omega},\ldots,\ma[N]{\Omega})+\ma{B}(\ve{x}),
\end{align*}
where $\ma{B}(\ve{x})$, the linearization of the homogeneous Lohe model \eqref{eq:homo}, is given in \cite{markdahl2018tac}. We can interpret the matrix $\ma{A}(\ve{x})$ as a the result of perturbing the matrix $\ma{B}(\ve{x})$.

Let $\beta(\ve{x})=\max_{\|\ve{v}\|_2=1}\langle\ve{v},\ma{B}(\ve{x})\ve{v}\rangle$ denote the spectral abscissa of $\ma{B}(\ve{x})$. By matrix perturbation theory \cite{kahan1975spectra},
\begin{align}
|\beta(\ve{x})-\re\alpha(\ve{x})|&\leq\|\ma{A}(\ve{x})-\ma{B}(\ve{x})\|_2\nonumber\\
&=\bigl(\sum_{i\in\V}\!\|\ma[i]{\Omega}\|^2_2\bigr)^{\frac12}.\label{eq:kahan}
\end{align}
A lower bound for $\beta(\ve{x})$ is given in \cite{markdahl2018tac}. Based on this bound we derive the following inequality which holds at any  dispersed equilibrium $\ve{x}$ (see Definition \ref{def:incohesive}),
\begin{align}\label{eq:result}
\beta(\ve{x})\geq\tfrac{K}{n+1}(n-1-\cos\tfrac{\pi}{N})(1-\cos\tfrac{\pi}N).
\end{align}
Suppose that the assumption \eqref{eq:th} of Theorem \ref{th:main} holds, then
\begin{align*}
|\beta(\ve{x})-\re\alpha(\ve{x})|&\leq\bigl(\sum_{i\in\V}\!\|\ma[i]{\Omega}\|^2_2\bigr)^{\frac12}\\
&<\tfrac{K}{n+1}(n-1-\cos\tfrac{\pi}{N})(1-\cos\tfrac{\pi}N)\\
&\leq\beta(\ve{x}),
\end{align*}
where we used \eqref{eq:kahan}, \eqref{eq:result}. From $|\beta(\ve{x})-\re\alpha(\ve{x})|<\beta(\ve{x})$ we get  $\re\alpha(\ve{x})>0$. It follows that the equilibrium $\ve{x}$ is exponentially unstable.\qedsymboltext

\begin{remark}
Theorem \ref{th:main} can be interpreted as a robustness result for a previous theorem on almost global synchronization of the homogeneous Lohe model \cite{markdahl2018tac}. We know that a set of perturbed equilibria that is close to the consensus manifold remains asymptotically stable \cite{lageman2016consensus}. In this paper we show that all other equilibria remain exponentially unstable for all $\ma[i]{\Omega}$ that satisfy \eqref{eq:th}. Compared to \cite{markdahl2018tac}, the property of almost global practical synchronization could not be established for the heterogeneous Lohe model. The result \cite{markdahl2018tac} utilize that  analytic gradient descent flows have very strong convergence properties \cite{lageman2007convergence}. These extend to some multiplicative perturbations, but not to additive perturbations like \eqref{eq:flow}. In general, there exist cases where a set of exponentially instable equilibria can be almost globally attractive \cite{freeman2013global}. If the dynamics belong to a class of nice functions (\eg gradients flows of analytic functions on compact sets), then such pathological behaviour is not possible. In our case, the right-hand side of \eqref{eq:dxidt} is not a gradient descent flow. In theory, it can not be ruled out that the behaviour of the heterogeneous Lohe model \eqref{eq:dxidt} is pathological. However, it is unlikely that \eqref{eq:dxidt} is a pathological case. Future work will explore the convergence behavoir of the Lohe model \eqref{eq:dxidt} near exponentially unstable equilibria.
\end{remark}

\section{Future work}

\noindent We provide an instability result based on the center manifold theorem for the heterogeneous Lohe model. The literature already contains local stability results for the desired equilibria, \ie practical synchronization for the perturbed consensus manifold \cite{lageman2016consensus}. Taken together, local stability of desired equilibria and local instability of undesired equilibria should amount to almost global asymptotical stability of the desired equilibria. However, some technical issues prevent this conclusion. The main question is if there exists trajectories that do not converge to equilibria, \eg limit cycles in the case of large $\ma[i]{\Omega}$, and how to account for them in the stability analysis. Future work will bridge this gap in our understanding and address the case of large perturbations.

\bibliographystyle{unsrt}
\bibliography{IEEEbib}

\appendix

\section{Proof of Theorem X}

\subsection{The center manifold theorem}

\label{sec:linearized}

\noindent The indirect method of Lyapunov states that a nonlinear system is unstable if the matrix that characterizes its linearization has an eigenvalue with strictly positive real part \cite{khalil2002nonlinear}. However, we must also account for the fact that \eqref{eq:dxidt} evolves on $\mathcal{S}^{n}\subset\R^{n+1}$. To this end, we use the center manifold theorem. Extend the heterogeneous Lohe model from $\mathcal{S}^n$ to an open neighborhood of $\mathcal{S}^n$ in $\R^{n+1}\backslash\{\ve{0}\}$ as
\begin{align}\label{eq:extension}
\vd[i]{v}&=\ma[i]{\Omega}\ve[i]{v}+(\ma[n]{I}-\tfrac{\ve[i]{v}}{\|\ve[i]{v}\|_2\!\!\!}\,(\tfrac{\ve[i]{v}}{\|\ve[i]{v}\|_2\!\!\!}\,\,)\mtr)\sum_{j\in\Ni}k_{ij}\tfrac{\ve[j]{v}}{\|\ve[j]{v}\|_2\!\!\!}\,\,,
\end{align}
where $\ve[i]{v}(0)\in\R^{n+1}$. Note that $\|\ve[i]{v}\|^2_2$ is constant along trajectories of the system since $\langle\vd[i]{v},\ve[i]{v}\rangle=0$. Also note that \eqref{eq:extension}  and the linearization of \eqref{eq:extension} on $(\mathcal{S}^n)^N$ equals the linearization of the heterogeneous Lohe model \eqref{eq:dxidt}.

The center manifold theorem states that the unstable manifold at an equilibrium $\ve{v}$ of \eqref{eq:extension} is tangent to the unstable subspace of the linearization of \eqref{eq:extension} \cite{sastry2013nonlinear}. Let $\ve{x}\in(\mathcal{S}^n)^N\subset\R^{N(n+1)}$ be an exponentially unstable equilibrium of the Lohe mode \eqref{eq:dxidt}, then $\ve{x}$ is also an exponentially unstable equilibrium of \eqref{eq:extension}. Part of the unstable manifold of \eqref{eq:extension} at $\ve{x}$ is tangent to $\ts[(\mathcal{S}^n)^N]{\ve{x}}$. Note that this part of the unstable manifold belongs to $(\mathcal{S}^n)^N$. To see this, follow a trajectory on the unstable manifold that satisfies $\lim_{t\rightarrow-\infty}\ve{v}(t)=\ve{x}$. Then  $\|\ve[i]{v}(t)\|_2=\|\ve[i]{x}\|_2=1$ or else we contradict $\langle\vd[i]{v},\ve[i]{v}\rangle=0$. As such, the equilibrium $\ve{x}$ of the Lohe model \eqref{eq:dxidt}, \ie the restriction of \eqref{eq:extension} to $(\mathcal{S}^n)^N$, is not Lyapunov stable.

\subsection{The linearized system}

\noindent Recall that the Lohe model \eqref{eq:dxidt} is the sum of a drift term and a gradient descent term, \ie $\vd[i]{x}=\ma[i]{\Omega}\ve[i]{x}-\grad_i V$. The map from a system to its linearization is itself linear, so we can linearize the two terms separetly. Note that $\vd[i]{y}=\ma[i]{\Omega}\ve[i]{y}$ is linear. The linearization of the gradient descent flow $\vd[i]{z}=-\grad_i V$ is given in \cite{markdahl2018tac}.

The linearization of $\vd[i]{z}=-\grad_i V$ is given by the $N(n+1)\times N(n+1)$ block matrix $\ma{B}(\ve{z})=[\ma[ij]{B}]$ where each $(n+1)\times(n+1)$ block is
\begin{align}\label{eq:B}
\ma[ij]{B}&=\begin{cases}
-\sum_{j\in\Ni}k_{ij}\langle\ve[j]{z},\ve[i]{z}\rangle(\ma[n]{I}-\ve[i]{z}\vet[i]{z})&\textrm{if }j=i,\\
k_{ij}(\ma[n]{I}-\ve[i]{z}\vet[i]{z})(\ma[n]{I}-\ve[j]{z}\vet[j]{z})&\textrm{otherwise.}
\end{cases}
\end{align}
The linearization of \eqref{eq:dxidt} is hence characterized by a block matrix $\ma{A}(\ve{x})=[\ma[ij]{A}]$ where the blocks are given by
\begin{align}\label{eq:A}
\ma[ij]{A}&=\begin{cases}
\ma[i]{\Omega}+\ma[ii]{B}&\textrm{if } j=i,\\
\ma[ij]{B} &\textrm{otherwise.}
\end{cases}
\end{align}

Equation \eqref{eq:A} allows us to see the linearization matrix $\ma{A}(\ve{x})$ as a perturbation of the matrix $\ma{B}(\ve{x})$,
\begin{align}\label{eq:perturbation}
\ma{A}(\ve{x})=\diag(\ma[1]{\Omega},\ldots,\ma[N]{\Omega})+\ma{B}(\ve{x}),
\end{align}
where we assume that $\|\diag(\ma[1]{\Omega},\ldots,\ma[N]{\Omega})\|_2$ is small. 
To analyze the stability of $\ma{A}(\ve{x})$, we first find a specific lower bound on the largest eigenvalue (the \emph{spectral abscissa})  of $\ma{B}(\ve{x})$. Then we use matrix perturbation theory in the case of small $\|\ma[i]{\Omega}\|_2$ to bound the spectral abscissa of $\ma{A}(\ve{x})$ below by a small positive number. The equilibrium $\ve{x}$ is then exponentially unstable by the direct method of Lyapunov.


\subsection{Eigenvalues of the symmetric matrix}

\noindent The spectral abscissa $\beta(\ve{x}):=\max_{\|\ve{v}\|_2=1}\langle\ve{v},\ma{B}(\ve{x})\ve{v}\rangle$ of $\ma{B}(\ve{x})$ given by \eqref{eq:B} is bounded below as
\begin{align*}
\beta(\ve{x})&\geq
\max_{\|\ve{w}\|_2=1}\tfrac{1}N\langle\ve{w},\sum_{i\in\V}\sum_{j\in\V}\ma[ij]{B}\ve{w}\rangle\\
&\geq\tfrac{1}{N(n+1)}\sum_{i\in\V}\sum_{j\in\V}\trace\ma[ij]{B}\geq\tfrac{1}{N(n+1)}\sum_{i\in\V}\sum_{j\in\Ni}\trace\ma[ij]{B}\\
&=\tfrac{2}{N(n+1)}\!\sum_{\{i,j\}\in\E}\! k_{ij}(-n\langle\ve[j]{x},\ve[i]{x}\rangle+n-1+\langle\ve[j]{x},\ve[i]{x}\rangle^2),
\end{align*}
where we have chosen $\ve{v}=\tfrac1N[\vet{w},\ldots,\vet{w}]$ for any $\ve{w}\in\R^{n+1}$ such that $\|\ve{w}\|_2=1$ whereby $\|\ve{v}\|_2=1$. This expression is strictly positive for all $n\geq2$ and all $(\ve[i]{z})_{i=1}^N\notin\mathcal{C}$. 
%
For future reference we define the lower bound of $\beta(\ve{x})$ as a function
\begin{align*}
f:\ve{x}\mapsto\tfrac{2}{N(n+1)}\!\sum_{\{i,j\}\in\E}\!k_{ij}(n-1-\langle\ve[j]{x},\ve[i]{x}\rangle)(1-\langle\ve[j]{x},\ve[i]{x}\rangle),
\end{align*}
\ie $\beta(\ve{x})\geq f(\ve{x})$ holds. Let $\theta_{ij}$ denote the angle between agent $i$ and $j$, \ie  $\cos\theta_{ij}=\langle\ve[j]{x},\ve[i]{x}\rangle$. Then we can write
\begin{align*}
f:\ve{x}\mapsto\tfrac{2}{N(n+1)}\!\sum_{\{i,j\}\in\E}\!k_{ij}(n-1-\cos\theta_{ij})(1-\cos\theta_{ij}).
\end{align*}

\subsection{Optimization of bound}

\noindent Note that $f$ depends on $\ve{x}$. To remove this dependence, we minimize $f$ over all dispersed equilibria. For an equilibrium to be regarded as dispersed, by Definition \ref{def:incohesive} we require that the set $\{\ve[1]{x},\ldots,\ve[N]{x}\}$ is not contained in an open hemisphere of $\mathcal{S}^n$. The minimization must be done both with respect to the graph $\mathcal{G}$ and the configuration $\ve{x}=(\ve[i]{x})_{i=1}^N$.

Each term in $f$ is positive, wherefore removing a link from $\mathcal{G}$ leads to a decrease in $f$. It follows that the graph which minimizes $f$ is a tree. By definition of $\ve{z}$ being dispersed, $\{\ve[1]{x},\ldots,\ve[N]{x}\}$ is not contained in a hemisphere. There must hence be a tuple $\mathcal{P}=(i_1,\ldots,i_k)$ for some $k\leq N$ such that $\sum_{j=1}^k\theta_{i_ji_{j+1}}\geq\pi$. Note that $f$ decreases with increasing $\theta_{ij}$. Hence, for any index $i\notin\mathcal{P}$, it is suboptimal to not set $\theta_{ik}=0$, where $k$ is any neighbor of $i$. There is no loss of optimality in assuming that $\mathcal{P}=(1,\ldots,N)$, \ie $\mathcal{G}$ is a path.


The choice of agent placements that minimize $f$ solves
\begin{align*}
\min f=\tfrac{2}{N(n+1)}&\sum_{i\in\V} k_{ij}(n-1-\cos\theta_{i,i+1})(1-\cos\theta_{i,i+1})\\
&\sum_{i\in\V}\theta_{i,i+1}=\pi,
\end{align*}
where addition of indices is modulo $N$. The inequality in the constraint has been replaced by an equality since $\mu$ decreases with increasing $\theta_{ij}$, \ie it is suboptimal to chose $(\ve[i]{x})_{i=1}^N$ such that $\sum_{i\in\V}\theta_{i,i+1}>\pi$. The coupling gains $k_{ij}$  complicate the next step of the analysis wherefore we let $K:=\min_{\{i,j\}\in\E}k_{ij}$ and replace the problem by
\begin{align}
\min g(\ve{x}):=\tfrac{2K}{N(n+1)}&\sum_{i\in\V} (n-1-\cos\theta_{i,i+1})(1-\cos\theta_{i,i+1})\nonumber\\
c(\ve{x}){}:=&\sum_{i\in\V}\theta_{i,i+1}-\pi=0,\tag{NP}\label{eq:NP}
\end{align}
where $\beta(\ve{x})\geq g(\ve{x})$ with equality if $k_{ij}=K$.

\subsection{Solution to optimization problem for large $N$}

\noindent To solve the nonlinear programming problem \eqref{eq:NP}, we use the Lagrange conditions for optimality \cite{nocedal1999numerical}. There is only one constraint $c(\ve{x})=0$, the gradient of which is non-zero. The optimal solution hence satisfies
\begin{align}\label{eq:lagrange}
\nabla g(\ve{x})+\lambda c(\ve{x})=\nabla g(\ve{x})+\lambda\ve{1}=\ve{0},
\end{align}
where $\lambda$ is a Lagrange multiplier and $\ve{1}=[1,\ldots,1]\mtr$. It follows that 
\begin{align}\label{eq:trig}
\sin\theta_{i,i+1}(n-2\cos\theta_{i,i+1})=-\tfrac{N(n+1)}{2a}\lambda,
\end{align}
where $\lambda\leq0$ is required for $\theta_{i,i+1}\in[0,\pi]$ to exist.

Both factors in the left-hand side of \eqref{eq:trig} increase on $[0,\pi/2]$ and decrease on $[\pi/2,\pi]$. As such, the curve of $\sin\theta_{i,i+1}(n-2\cos\theta_{i,i+1})$ can intersect a constant at most twice on $[0,\pi]$. Moreover, one of the solutions is larger than $\pi/2$. However, since $\sum_{i\in\V}\theta_{i,i+1}=\pi$, at most one $\theta_{i,i+1}$ is larger than $\pi/2$. There are hence only two solutions to \eqref{eq:lagrange}: either $\theta_{i,i+1}=\pi/N$ for all $i\in\V$ or $\theta_{j,j+1}=\varphi>\pi/2$, $\theta_{i,i+1}=(\pi-\varphi)/(N-1)$, for all $i\in\V\backslash\{j\}$ (this is one solution up to permutations of $\{1,\ldots,N\}$). 

Consider the case where $\theta_{i,i+1}=\pi/N$. Then the objective function value is 
\begin{align*}
g_1&:=\tfrac{2K}{n+1}(n-1-\cos\tfrac{\pi}N)(1-\cos\tfrac{\pi}N).
\end{align*}
In the case when $\theta_{j,j+1}=\varphi$, then the objective function value is
\begin{align*}
g_2:={}&\tfrac{2K(N-1)}{N(n+1)}(n-1-\cos\tfrac{\pi-\varphi}{N-1})(1-\cos\tfrac{\pi-\varphi}{N-1})+\\
&\tfrac{2K}{N(n+1)}(n-1-\cos\varphi)(1-\cos\varphi)
\end{align*}
The value $g_1$ belongs to $\mathcal{O}(K/N^2)$ since  $1-\cos\pi/N\in\mathcal{O}(1/N^2)$ if $n\geq3$ (or $\mathcal{O}(K/N^4)$ if $n=2$) whereas $g_2$ belongs to $\mathcal{O}(K/N)$ due to the second term being bounded below by $2K(n-1)/(N(n+1))$. As such, $g_1$ is optimal for sufficiently large $N$. Note that $\theta_{i,i+1}=\pi/N<1$ for $N\geq5$ and
\begin{align*}
g_1&\leq\tfrac{2K}{n+1}(n-1)\tfrac12\tfrac{\pi^2}{N^2}<\tfrac{2K}{n+1}(n-1)\tfrac{1}N<g_2
\end{align*}
due to $\pi^2\leq2N$. Hence $\theta_{i,i+1}=\pi/N$ is optimal for $N\geq5$. 

\subsection{Three special cases}

\noindent Consider the remaining cases of $N\in\{2,3,4\}$. If $N=2$, then
\begin{align*}
g_1={}&\tfrac{2K}{n+1}(n-1),\\
g_2={}&\tfrac{K}{n+1}(n-1+\cos\varphi)(1+\cos\varphi)+\\
&\tfrac{K}{n+1}(n-1-\cos\varphi)(1-\cos\varphi)\\
={}&g_1+\tfrac{2K}{n+1}\cos^2\varphi.
\end{align*}
If $N=3$, then
\begin{align*}
g_1={}&\tfrac{K}{n+1}(n-\tfrac32),\\
g_2={}&\tfrac{4K}{3(n+1)}(n-1-\sin\tfrac{\varphi}2)(1-\sin\tfrac{\varphi}2)+\\
&\tfrac{2K}{3(n+1)}(n-1-\cos\varphi)(1-\cos\varphi)\\
={}&\tfrac{K}{n+1}[n(2-\tfrac{2}{3}(2\sin\tfrac{\varphi}2+\cos\varphi))\\
&\phantom{\tfrac{a}{n+1}(}-2(1-\tfrac{1}{3}(1-\cos\varphi+\cos^2\varphi))]\\
>{}&\tfrac{K}{n+1}(n-\tfrac43)>g_1,
\end{align*}
where we omitted a few steps. Finally, if $N=4$, then
\begin{align*}
g_1={}&\tfrac{2K}{n+1}(n-1-\tfrac{1}{\sqrt2})(1-\tfrac{1}{\sqrt2})\\
={}&\tfrac{K}{n+1}((2-\sqrt2)n-1),\\
g_2={}&\tfrac{3K}{2(n+1)}(n-1-\cos\tfrac{\pi-\varphi}3)(1-\cos\tfrac{\pi-\varphi}3)+\\
&\tfrac{K}{2(n+1)}(n-1-\cos\varphi)(1-\cos\varphi)\\
={}&\tfrac{K}{n+1}[(2-\tfrac54\cos\varphi-\tfrac{3\sqrt{3}}4\sin\varphi)n-\\
&\phantom{\tfrac{K}{n+1}\,\,}(\tfrac32\sin^2\tfrac{\pi-\varphi}3+\tfrac12\sin^2\varphi)]\\
>{}&\tfrac{K}{n+1}((2-\tfrac{3\sqrt{3}}4)n-\tfrac12(3\sin^2\tfrac{\pi-\varphi}3+\sin^2\varphi))\\
>{}&\tfrac{K}{n+1}((2-\tfrac{3\sqrt{3}}4)n-\tfrac78)>g_1,
\end{align*}
where we have omitted a few steps and $\sqrt{2}>3\sqrt{3}/4$ can be verified by numerical calculation.

\subsection{Matrix perturbation theory}

\noindent We have a lower bound for the spectral abscissa of the linearization of the homogeneous Lohe model at any dispersed equilibria. It remains to relate this result to the inhomogeneous Lohe model. This is done via the characterization \eqref{eq:perturbation} of the linearization of the inhomogenous model as a perturbation of the homogeneous model. 

A \emph{Hermitian matrix} $\ma{X}$ equals its conjugate transpose, \ie $\ma{X}^*=\ma{X}$. Given an Hermitian matrix $\ma{X}$ and an arbitrary matrix $\ma{Y}$, the following result from matrix perturbation theory relates the spectrum of $\ma{X}+\ma{Y}$ to the spectrum of $\ma{X}$ and $\ma{Y}$:


\begin{theorem}[Kahan \cite{kahan1975spectra}]\label{th:Kahan} Let $\ma{X}\in\C^{n\times n}$ be a Hermitian matrix with eigenvalues $\chi_1\leq\ldots\leq\chi_n$ and $\ma{Y}\in\C^{n\times n}$ be an arbitrary matrix. Let $\sigma_j$ with $\re\sigma_1\leq\ldots\leq\re\sigma_n$ denote the eigenvalues of $\ma{X}+\ma{Y}$. Then
\begin{align*}
(\sum_{j=1}^n(\chi_j-\re\sigma_j)^2)^\frac12\leq{}&\tfrac12\|\ma{Y}+\ma{Y}^*\|_2+\\
&(\tfrac14\|\ma{Y}-\ma{Y}^*\|_2^2-\sum_{j=1}^n(\im\sigma_j)^2)^\frac12.
\end{align*}
\end{theorem}

Note that for the case we are interested in, $\ma{Y}\in\so$, \ie $\mat{Y}=-\ma{Y}$. The inequality in Theorem \ref{th:Kahan} allows several simplifications for $\ma{Y}\in\so$, including
\begin{align*}
|\chi_n-\re\sigma_n|\leq\|\ma{Y}\|_2.
\end{align*}
Note that a sharper inequality could be obtainable if we knew more about $\chi_j$ and $\sigma_j$. However, we are satisfied with this rather elegant bound.

%
%
%

The spectral abscissa $\beta(\ve{x})$ of the linearization matrix $\ma{B}(\ve{x})$ of the homogeneous Lohe model given by \eqref{eq:B} satisfies 
\begin{align}\label{eq:step1}
\beta(\ve{x})\geq \tfrac{2K}{n+1}(n-1-\cos\tfrac{\pi}N)(1-\cos\tfrac{\pi}N).
\end{align}
The linearization matrix $\ma{A}(\ve{x})$ of the heterogenous Lohe model is separated from $\ma{B}(\ve{x})$ by
\begin{align*}
\ma{A}(\ve{x})-\ma{B}(\ve{x})=\diag(\ma[1]{\Omega},\ldots,\ma[N]{\Omega}).
\end{align*}
Let $\re\alpha(\ve{x})$ denote the spectral abscissa of $\ma{A}(\ve{x})$. By the matrix perturbation result, Theorem \ref{th:Kahan}, we have
\begin{align}\label{eq:step2}
|\beta(\ve{x})-\re\alpha(\ve{x})|&\leq \|\diag(\ma[1]{\Omega},\ldots,\ma[N]{\Omega})\|_2.
\end{align}
Note that $\|\diag(\ma[1]{\Omega},\ldots,\ma[N]{\Omega})\|_2=(\sum_{j\in\V}\|\ma[j]{\Omega}\|_2^2)^\frac12$. By the assumption of Theorem \ref{th:main},
\begin{align}\label{eq:step3}
(\sum_{j\in\V}\|\ma[j]{\Omega}\|_2^2)^\frac12<\tfrac{2K}{n+1}(n-1-\cos\tfrac{\pi}N)(1-\cos\tfrac{\pi}N).
\end{align}
Combine \eqref{eq:step1}--\eqref{eq:step3} to find $|\beta(\ve{x})-\re\alpha(\ve{x})|<\beta(\ve{x})$, which implies $\re\alpha(\ve{x})>0$. Since $\alpha(\ve{x})$ is an eigenvalue of the linearization matrix $\ma{A}(\ve{x})$ of the heterogeneous Lohe model at any dispersed equilibrium $\ve{x}$, it follows that the heterogeneous Lohe model is exponentially unstable at $\ve{x}$.

\end{document}